\documentclass[10pt,a4paper,reqno]{amsart}
\usepackage{amsfonts,amsthm,latexsym,amsmath,amssymb,amscd,amsmath,epsf, 
graphicx}

\newtheorem{tm}{Theorem}

\newtheorem{rem}[tm]{Remark}

\newtheorem{lm}[tm]{Lemma}

\newtheorem{cor}[tm]{Corollary}
\newtheorem{prop}[tm]{Proposition}
\newtheorem{nota}[tm]{Notation}

\newcommand{\la}{\lambda}

\newcommand{\C}{\mathcal C}

\begin{document}
\title{Narayana numbers and Schur-Szeg\"o composition}

\author[V.~Kostov]{Vladimir P. Kostov }
\address{Universit\'e de Nice, Laboratoire de Math\'ematiques, 
Parc Valrose, 06108 Nice Cedex 2, France}
             \email{kostov\@math.unice.fr}

\author[B.~Shapiro]{Boris Z. Shapiro}
\address{Department of Mathematics, Stockholm University, SE-106 91
Stockholm,
         Sweden}
\email{shapiro\@math.su.se}

\begin{abstract}
In the present paper we find a new interpretation of Narayana polynomials 
$N_n(x)$ which are the generating polynomials for the Narayana numbers 
$N_{n,k}$ counting Dyck paths of length $n$ and with exactly $k$ peaks, 
see e.g. \cite {Su}. 
(These numbers appeared recently in a number of different combinatorial 
situations, \cite {DSV,STT,YY}.) Strangely enough Narayana polynomials 
also occur 
as  limits as $n\to \infty$ of the sequences of eigenpolynomials of  
the Schur-Szeg\"o composition map sending $(n-1)$-tuples 
of polynomials of the form $(x+1)^{n-1}(x+a)$ to their Schur-Szeg\"o product, 
see below. As a corollary we obtain that every $N_n(x)$ has all roots real and 
non-positive. Additionally, we present an explicit formula for the 
density and the distribution function of the asymptotic root-counting 
measure of the polynomial sequence $\{N_n(x)\}$.

\end{abstract}



\date{}
\maketitle 

\section{Introduction}

\subsection{The Narayana numbers, triangle and  polynomials}

The {\em Narayana numbers} $N_{n,k},\; 1\le k\le n,$ apparently 
introduced by G.~Kreweras in \cite {Kr} are given by: 
$$N_{n,k}=\frac{1}{n}C_n^{k-1}C_n^k,$$
where $C_j^i$ stands for the usual binomial coefficient, i.e. 
$C_j^i=\frac{j!}{i!(j-i)!}$. 

The latter formula immediately implies that for any fixed $k$ the Narayana 
numbers $N_{n,k}$ are given by a polynomial in $n$ of degree $2k-2$ divisible 
by $n$. It is known that $N_{n,k}$ counts, in particular, 
the number of expressions containing $n$ pairs of parentheses  which are 
correctly matched and which contain exactly $k$ distinct nestings and also 
the number of Dyck paths of length $n$ with exactly $k$ peaks. (Recall that 
a Dyck path is a staircase walk from $(0,0)$ to $(n,n)$ that lies strictly 
above (but may touch) the diagonal $y=x$.) Some other combinatorial 
interpretations of $N_{n,k}$ can be found in \cite{Su} and references  
therein. 

The triangle 
\begin{equation}\label{NT} 
\begin{array}{cccccccccc}&&&&1&&&&\\&&&1&&1&&&\\
&&1&&3&&1&&\\
&1&&6&&6&&1&\\1&&10&&20&&10&&1&
\end{array}
\end{equation} 

\medskip
\noindent
of Narayana numbers $N_{n,k}$ read by rows is called the {\em Narayana  
triangle}. (Later we will also interpret this triangle as an infinite 
lower-triangular matrix taking its left side of ones as the first column and 
its right side of ones  as the main diagonal, see (\ref{NNT}).)  

The generating functions of the rows of the above triangle are called the 
{\em Narayana polynomials}. More exactly, following the standard convention  
(see \cite{B})  one defines  the $n$-th {\em Narayana polynomial} 
by the formula

\[ N_n(x)=\sum _{k=1}^{n}N_{n,k}x^k~~.\] 

In what follows we will use the following notions. If $P(x)$ is a 
univariate polynomial of degree $n$, then its {\em reversion} or the 
{\em reverted polynomial} $P^R(x)$ is defined as $P^R(x)=x^nP(1/x)$. 
A polynomial $P(x)$ is called  {\em self-reciprocal} if it coincides with  
its revertion up to a sign, i.e. 
$P(x)=\pm P^R(x)$. Hence for any self-reciprocal $P(x)$ if $P(x)$ vanishes at 
$x_0$ then $P(x)$ vanishes at $1/x_0$ as well. A polynomial $P(x)$ is called 
{\em hyperbolic} if all its roots are real. 

\begin{rem}\label{Narrem}
Each polynomial $N_n(x)$ has a simple root at $0$ and  each  $N_n(x)/x$ is 
self-reciprocal.  
\end{rem}

The following simple $3$-term recurrence relation satisfied by Narayana 
polynomials was found in  \cite[p. 2]{Su}:

\begin{equation}\label{recur} 
(n+1)N_n(x)=(2n-1)(1+x)N_{n-1}(x)-(n-2)(x-1)^2N_{n-2}(x), 
\end{equation} 
with the initial conditions $N_1(x)=x,\; N_2(x)=x^2+x$. 



\subsection{Schur-Szeg\"o composition}

The {\em  Schur-Szeg\"o composition  (CSS)} of two degree $n$ polynomials 
$P=\sum _{j=0}^np_jx^j$ and $Q=\sum _{j=0}^nq_jx^j$ is defined by the formula: 
$$P~^{\ast}_{n}~Q=\sum _{j=0}^np_jq_jx^j/C_n^j.$$
When the same $P$ and $Q$ are considered as polynomials of degree $n+k$  
with vanishing $k$ 
leading coefficients  then in accordance with the above formula one gets: 
$$P~^{~\, \ast}_{n+k}~Q=\sum _{j=0}^np_jq_{j}x^j/C_{n+k}^j.$$
 Extending these formulas one defines the composition of $s$ polynomials 
by the formula: 

\[ P_1~^{~\, \ast}_{n+k}~\cdots ~^{~\, \ast}_{n+k}~P_s=
\sum _{j=0}^np_{1,j}\cdots p_{s,j}x^j/(C_{n+k}^j)^{s-1}~.\] 
(For more details on CSS see \cite{Pr,RS}.)  

Our main goal below will be a further study of a certain linear 
inhomogeneous  map $\Phi_n$ initially considered in \cite{Ko2,AlKo}.  
Namely, in these papers the first author of the present paper has shown  
the possibility to present every monic polynomial of degree $n$ with complex 
coefficients and  vanishing  at $(-1)$ in the form:   

\begin{equation}\label{fact}
P=K_{a_1}~^{\ast}_{n}~\cdots ~^{\ast}_{n}~K_{a_{n-1}}
\end{equation} 
where each {\em composition factor} 
$K_{a_i}$ equals $(x+1)^{n-1}(x+a_i),\; a_i\in \mathbf C.$ 
(For the sake of convenience, we set $K_{\infty}:=(x+1)^{n-1}$.) 
Now we can introduce the map $\Phi_n$. 

\begin{nota}\label{sigmanu}
For any  $P(x):=(x+1)(x^{n-1}+c_1x^{n-2}+\cdots +c_{n-2}x+c_{n-1})$ and for 
$\nu =1,\ldots ,n-1$ set  
$\sigma _{\nu}:=\sum _{1\leq j_1<\cdots <j_{\nu}\leq n-1}
a_{j_1}\cdots a_{j_{\nu}}$, i.e.  define $\sigma _{\nu}$ as the $\nu$-th 
elementary symmetric function of the roots of the composition 
factors presenting $P(x)$. Finally,  denote by $\Phi_n$ 
the mapping $(c_1,\ldots ,c_{n-1})\mapsto (\sigma _1,\ldots ,\sigma _{n-1})$. 
\end{nota}

Obviously, $\Phi_n$ is   linear inhomogeneous.   The following theorem was 
proven in  \cite{Ko1}.

\begin{tm}\label{maintm}
\begin{enumerate}
\item The mapping $\Phi_n$ has $n-1$ distinct real eigenvalues  
$\lambda _{1,n}=1$, $\lambda _{2,n}=\frac{n}{n-1}$, 
$\lambda _{3,n}=\frac{n^2}{(n-1)(n-2)}$, 
$\ldots$, $\lambda _{n-1,n}=\frac{n^{n-2}}{(n-1)!}$. 
\noindent
\item The corresponding eigenvectors are  monic  polynomials of degree 
$n-1$ vanishing at $(-1)$ and have the form:   
$(x+1)^{n-1}$, $x(x+1)^{n-2}$, $x(x+1)^{n-3}Q_{1,n}(x)$, 
$\ldots$, $x(x+1)Q_{n-3,n}(x)$ where 
$\deg Q_{j,n}(x)=j$, $j=1,\ldots ,n-3$, $Q_{j,n}(-1)\neq 0$.  
The coefficients of each polynomial $Q_{j,n}(x)$ are rational numbers. 
\noindent
\item Each  $Q_{j,n}(x)$ is self-reciprocal. More exactly, 
$(Q_{j,n}(x))^R=(-1)^jQ_{j,n}(x)$.
\noindent
\item  The roots of each $Q_{j,n}(x)$, $1\leq j\leq n-3$, are positive and 
distinct.
\noindent
\item  For $j$ odd (resp. for $j$ even) one has  $Q_{j,n}(1)=0$ (resp. 
$Q_{j,n}(1)\neq 0$). 
Additionally, the middle coefficient in $(x+1)^{n-j-2}Q_{j,n}(x)$ vanishes  
if $n$ is even and $j$ is odd. 
\noindent
\item For any $j$ fixed and $n\rightarrow \infty$ the sequence of polynomials 
$Q_{j,n}(x)$ converges coefficientwise to the monic  polynomial $Q_{j}^*(x)$ 
of degree $j$ which has  rational coefficients, all roots positive, and 
satisfies the equality $(Q_{j}^*(x))^R=(-1)^jQ_{j}^*(x)$ 
and the condition $Q_{j}^*(1)=0$ for $j$ odd.
\end{enumerate}
\end{tm}

\begin{rem} In Theorem~\ref{maintm} we consider the action of $\Phi_n$ of the 
affine $(n-1)$-dimensional  space of all monic polynomials of degree $n-1$. If 
we extend this action to the ambient linear space of all polynomials of degree 
at most $(n-1)$ then we acquire one more eigenvalue and eigenvector. Namely, 
the polynomial $(x+1)^{n-2}$ is the eigenvector of $\Phi_n$ with the 
eigenvalue $1$.
\end{rem}

\subsection {Main results} 

Set $M_j(x)=(-1)^{j-1}xQ_{j-1}^*(-x)$. The most important  
result of the present paper (with a rather lengthy proof)  is as follows, see 
details in Subsection~\ref{prMjNj}. 

\noindent
\begin{tm}\label{MjNj} For any positive integer $j$ the polynomial $M_j(x)$ 
coincides with the Narayana polynomial $N_j(x)$.
\end{tm}

Part (6) of the above Theorem~\ref{maintm} then  implies the following. 

\begin{cor}\label{hyperbol}
Narayana polynomials are hyperbolic for any $n\geq 1$.
\end{cor}

More information about the roots of $N_n(x)$ is given below. For its proof 
consult Subsection~\ref{printerlace}. 

\noindent 
\begin{tm}\label{interlace}
\begin{enumerate} 
\item   The number $(-1)$ is a simple root of $N_n(x)$ for any positive even 
integer $n$. For $n$ odd one has  $N_n(-1)\neq 0$;
\noindent
\item All roots of $N_n(x)$ are distinct and nonpositive; 
\noindent
\item The roots of $N_{n-1}(x)/x$ interlace with the ones   
of $N_{n}(x)/x$. Except for the origin the polynomials 
$N_{n-1}(x)$ and $N_n(x)$ have no root in common. 
\end{enumerate}
\end{tm}

Our final result is as follows, see details in Subsection~\ref{asymptotics}.  
Given a polynomial $P(x)$ of degree $l$ define its {\em root-counting measure} 
$\mu_P=\frac{1}{l}\sum_{i=1}^l\delta(x-x_i)$ where $\{x_1,\ldots ,x_l\}$ 
is the set of all roots of $P(x)$ listed with possible repetitions (equal to 
the respective multiplicities) and $\delta(x-x_i)$ is the standard Dirac 
delta-function supported at $x_i$. Given a sequence $\{P_n(x)\}, \deg P_n(x)=n,
\; n=1,2,\ldots$ we call  {\em asymptotic root-counting measure} of this 
sequence the weak limit $\mu=\lim_{n\to \infty} \mu_{P_n}$ (if it exists) 
understood in the sense of distribution theory. 

\begin{figure}[!htb]
\centerline{\hbox{\epsfysize=3.5cm\epsfbox{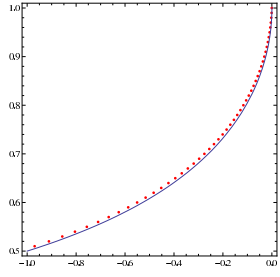}}}
\caption{The theoretical distribution $\kappa(x)$ and the empirical 
distribution of roots of $N_{100}(x)$ on the interval $[-1,0]$.}
\label{fig1}
\end{figure}

\begin{tm}\label{dens} The density $\rho(x)$ and the distribution function 
$\kappa(x)$ of the asymptotic 
root-counting measure of the sequence $\{N_n(x)\}$ of the Narayana polynomials 
are given by:
\begin{equation}\label{eq:dens}
\rho(x)=\frac{1}{\pi}\frac{1}{(1-x)\sqrt{-x}};
\quad\kappa(x)=1-\frac{2}{\pi}\arctan{\sqrt{-x}},\;x\le 0.
\end{equation}
\end{tm}

\begin{rem} Notice that self-reciprocity of $N_n(x)$ translates in the 
following (easily testable) property of $\rho(x)$:
$$x^2\rho(x)=\rho\left(\frac{1}{x}\right).$$ 

\end{rem}

\noindent
{\bf Acknowledments.} The authors are grateful to Professor Andrei 
Martinez-Finkel\-shtein for important conversations.  Research of the first author was partially supported by project 20682
of cooperation 
between CNRS and FAPESP "Zeros of algebraic polynomials. 
The second author expresses his  gratitude to  
Laboratoire de Math\'ematiques, Universit\'e de Nice-Sophia-Antipolis for the 
financial support of his visit to Nice in January 2008.  

\section{Proofs}

\subsection{Preliminaries}

To prove Theorems~\ref{MjNj} and \ref{interlace} we will need a detailed study 
of the map $\Phi_n$ and,  especially, of the equations defining its 
eigenvectors which in their turn give our polynomials $Q_{j,n}(x)$.   
\begin{nota}\label{psikj}
Set  $e _k(j)=\sigma_k(1,2,\ldots ,j)$ to be the value of the 
$k$-th  symmetric function on the $j$-tuple of numbers $(1, 2, \ldots, j)$, 
i.e. $e _k(j)=\sum _{1\leq \nu _1<\cdots <\nu _k\leq j}\nu _1\cdots \nu _k$, 
$k=1$, $\ldots$, $j$. Denote by $\phi _k(j)$ the sum $1^k+2^k+\cdots +j^k$. 
\end{nota}

\begin{rem}\label{divisible}
The quantity $e _k(j)$ (resp. $\phi _k(j)$) is a polynomial in $j$ of degree 
$2k$ (resp. of degree $k+1$) divisible by $j(j+1)$.
\end{rem}




Let $Q_{j,n}(x):=x^j+q_1x^{j-1}+\cdots +q_{j-1}x+(-1)^j$ be the polynomial 
introduced in Theorem~\ref{maintm} and set $l_j=(n-1)\cdots (n-j)$. 
Then by (1) of  Theorem~\ref{maintm} one has 
$$\lambda _{j+2,n}=n^{j+1}/l_{j+1}.$$ 
(the coefficients $q_{\nu}$ depend also on $j$ and $n$, but we prefer to 
avoid  double indices.)  
By  definition the polynomial $Q_{j,n}(x)$ satisfies the following relation:

\[ x(x+1)^{n-j-2}Q_{j,n}(x)=\] 
\[=\lambda _{j+2,n}
x(x+1)^{n-1}~^{\ast}_{n}~(x+a_1)(x+1)^{n-1}~^{\ast}_{n}~\cdots ~^{\ast}_{n}~
(x+a_j)(x+1)^{n-1}~^{\ast}_{n}~(x+1)^{n-1},\] 
where $\{-a_1,\ldots ,-a_j\}$ is the set of all roots of $Q_{j,n}(x)$. 
After  multiplication of both sides  of the latter relation by $l_{j+1}$ one 
gets that the coefficient $R_k$ of  $x^k$, $k\geq 1$, 
in the right-hand side equals 

\[ R_k:=n^{j+1}C_{n-1}^{k-1}(C_{n-1}^ka_1+C_{n-1}^{k-1})\cdots 
(C_{n-1}^ka_j+C_{n-1}^{k-1})C_{n-1}^k/(C_n^k)^{j+1}~.\] 
The corresponding coefficient  $L_k$ in the left-hand side equals 

\begin{equation}\label{Lk} 
L_k:=(n-1)\cdots (n-j-1)((-1)^jC_{n-j-2}^{k-1}+C_{n-j-2}^{k-2}q_{j-1}
+\cdots +C_{n-j-2}^0q_{j-k+1}).
\end{equation}

Therefore one has $q_{\nu}=\sigma _{\nu}$ 
(see Notation~\ref{sigmanu}) and, finally,

\begin{equation}\label{Rk} 
R_k=n^{j+1}C_{n-1}^{k-1}C_{n-1}^k((C_{n-1}^{k-1})^j+
\sum _{\nu =1}^{j-1}(C_{n-1}^{k-1})^{j-\nu}(C_{n-1}^k)^{\nu}q_{\nu}+
(-1)^j(C_{n-1}^k)^{j})/(C_n^k)^{j+1}.
\end{equation} 
Thus the coefficients $q_{\nu},\;\nu=1,\ldots ,j-1$ of $Q_{j,n}(x)$ 
solve the system of equations 

\begin{equation}\label{systSigma} 
(\Sigma )~~:~~\{ ~L_k=R_k~,~k=1,\ldots ,n-1~\}. 
\end{equation} 

\begin{lm}\label{seriesn-1}
The coefficients $q_{\nu}$ can be expanded in convergent series: 

\begin{equation}\label{qTaylor} 
q_{\nu}=q_{\nu}^{(0)}+\frac{q_{\nu}^{(1)}}{n-1}+\frac{q_{\nu}^{(2)}}{(n-1)^2}+
\cdots
\end{equation}
 with respect to  $\frac{1}{n-1}$, where the numbers 
$q_{\nu}^{(i)}\in {\bf R}$ are uniquely defined and independent of $n$. 
\end{lm}

\noindent
{\em Proof:} The coefficients $q_{\nu}$ solve  system $(\Sigma )$. They are 
uniquely defined because the polynomials $Q_{j,n}(x)$ are uniquely defined by 
the eigenvectors of the mapping $\Phi_n$. They can be expanded in 
convergent series 
in $\frac{1}{n-1}$ because the same property holds for the coefficients of 
system $(\Sigma )$. 
As $q_{\nu}$ are uniquely defined, thus $q_{\nu}^{(i)}$ are also uniquely 
defined.~~~~~$\Box$ 

\begin{rem}
We choose to expand $q_{\nu}$ as a series in $\frac{1}{n-1}$ (and not in 
$\frac{1}{n}$)  
because the eigenpolynomials of $\Phi_n$ (see (2) of Theorem~\ref{maintm}) 
are all of degree 
$n-1$. Besides, numerical computations show that it is the factor $n-1$ and 
not $n$ which appears most often in the denominators of the eigenvectors of 
the mapping $\Phi_n$.
\end{rem}

\begin{prop}\label{qj-1}
One has $q_1^{(0)}=(-1)^jq_{j-1}^{(0)}=-j(j+1)/2$. 
\end{prop}

\noindent
{\em Proof:} For $k=1$ one has 
\[ L_1=(-1)^j(n-1)\cdots (n-j-1)~~,~~ 
R_1=(n-1)
(1+\sum _{\nu =1}^{j-1}(n-1)^{\nu}q_{\nu}+(-1)^j(n-1)^j)~.\] 
The equality $L_1=R_1$ can be written in the form 
\[ (-1)^j(n-1)\cdots (n-j-1)=(-1)^j(n-1)^{j+1}+(n-1)^{j}q_{j-1}^{(0)}+
o((n-1)^{j})~.\] 
Hence $q_{j-1}^{(0)}+o(1)=(-1)^j((n-1)\cdots (n-j-1)-(n-1)^{j+1})/(n-1)^{j}$. 
Observe that 
\[ (n-1)\cdots (n-j-1)=(n-1)^{j+1}-(1+2+\cdots +j)(n-1)^j+o((n-1)^j)~.\]  
The quantity $q_{j-1}^{(0)}$ depends on $j$, but not on $n$. 
Therefore $(-1)^jq_{j-1}^{(0)}=-(1+2+\cdots +j)=-j(j+1)/2$.~~~~~$\Box$ 

\medskip 
The next statement is central. 

\begin{prop}\label{polynomial}
\begin{enumerate}
\item For each $(\nu ,i)$ fixed the coefficient  $q_{\nu}^{(i)}$  is given by  
a real polynomial in  $j$ of degree $2(\nu +i)$. 
\item For $i=0$ this 
polynomial is divisible by $j(j+1)$.
\end{enumerate} 
\end{prop}

\subsection{Proof of Proposition~\protect\ref{polynomial}
\protect\label{prpolynomial}}

$1^0$. To prove part (1) of the proposition we use induction on $\nu +i$. 
Proposition~\ref{qj-1}  constitutes  
the base of induction. The  step of induction is explained in $2^0$ -- $3^0$. 

Recall that the coefficients $q_{\nu}$ give the unique solution to  system 
$(\Sigma )$. From now on we assume that system 
$(\Sigma )$ is infinite, i.e. $k=1,2,\ldots$. Substituting  
  the expansions (\ref{qTaylor}) of the coefficients  $q_{\nu}$ in 
$(\Sigma )$ we obtain a new 
system (denoted by $(x)$) with  variables 
$q_{\nu}^{(i)}$, $\nu =1,\ldots ,j-1$, $i=0,1,\ldots$. After this substitution 
 the equation $L_k=R_k$ of system $(\Sigma )$ transforms into  an 
equation of the form $\sum _{l=l_0}^{\infty}A_{k,l}/(n-1)^l=0$ where the 
quantities $A_{k,l}$ are some linear inhomogeneous  functions of the variables 
$q_{\nu}^{(i)}$.  
(Notice that $A_{k,l}$ depend on $j$ but not on $n$.)
The latter equation holds for all $n\in {\bf N}$ if and only if all 
 $A_{k,l}$ vanish. (The equation $\{A_{k,l}=0\}$ is denoted by $(A_{k,l})$.) 

$2^0$. The solution to 
system $(x)$ is unique (which follows from the uniqueness of the polynomials 
$Q_{j,n}(x)$ for every fixed $n$, see Theorem~\ref{maintm}). 
This solution depends only on $j$. The 
self-reciprocity of $Q_{j,n}(x)$ implies that $q_{\nu}^{(i)}=
(-1)^jq_{j-\nu}^{(i)}$.

In what follows we consider subsystems of system $(x)$ of the form 
$\{ (A_{k,l}), l=l_0,\ldots ,l_1\}$, i.e. systems defined in accordance with 
the filtration of the space of Laurent series in $\frac{1}{n-1}$ by the degree 
of $\frac{1}{n-1}$. We set $l=s-j-k$.

\begin{nota}
Denote by ${\mathcal I}_{a,b}$ 
the set of  variables $\{q_a^{(0)}$, $q_{a+1}^{(1)}$, $\ldots$, 
$q_{a+b}^{(b)}\}$. 
\end{nota}

To settle part (1) of Proposition~\ref{polynomial} we need Lemmas~\ref{A} and 
\ref{poldeg} whose proofs are given after that of Proposition~\ref{polynomial}.
 
\begin{lm}\label{A}
The linear inhomogeneous form $A_{k,s-j-k}$ depends only on the variables in 
the set  ${\mathcal J}_{j-s,s-1}:={\mathcal I}_{j-s,s-1} 
\cup {\mathcal I}_{j-s+1,s-2}\cup \cdots \cup {\mathcal I}_{j-1,0}.$
\end{lm}


$3^0$. Suppose that the variables belonging to the set 
${\mathcal J}_{j-s+1,s-2}$ 
are already determined. (For $s=2$ one has 
${\mathcal J}_{j-s+1,s-2}={\mathcal I}_{j-s+1,s-2}={\mathcal I}_{j-1,0}=
\{ q_{j-1}^{(0)}\}$; 
see Proposition~\ref{qj-1}.) The system of $s$ linear equations 
$(A):=\{ (A_{k,s-j-k}), k=1,\ldots ,s\}$ is a system 
with $s$ unknown variables, namely, those in the set  
${\mathcal I}_{j-s,s-1}$. This 
system has a unique solution (which follows from the existence and uniqueness 
of the polynomials $Q_{j,n}(x)$, see Theorem~\ref{maintm}). 
Hence the variables in the set ${\mathcal I}_{j-s,s-1}$ are uniquely defined. 
  
\begin{lm}\label{poldeg}
The solution to system $(A)$ is an $s$-vector  consisting of 
real polynomials in $j$  of degree $2s$.
\end{lm}

This concludes the proof of the step of induction in part (1) of 
Proposition~\ref{polynomial}. 

$4^0$. For $\nu =1$ part (2) of the proposition follows from 
Proposition~\ref{qj-1}. When solving the linear system $(x)$ 
we express the variables in the set ${\mathcal I}_{j-s,s-1}$ as affine 
functions of the ones in the set 
${\mathcal J}_{j-s+1,s-2}$. Suppose that all variables in that set are shown 
to be polynomials divisible by $j(j+1)$. Then the variables in the set 
${\mathcal I}_{j-s,s-1}$ will be divisible by $j(j+1)$ if and only if this 
is the case of the constant terms of system $(A)$ 
(we call them CTs for short). 

The CTs are the coefficients of $(n-1)^{s-j-k}$ of the expression 
$(-1)^jU_1-U_2-(-1)^jU_3$ where 

\[ U_1=(n-1)\cdots (n-j-1)C_{n-j-2}^{k-1}~~,~~
U_2=n^{j+1}C_{n-1}^{k-1}C_{n-1}^k(C_{n-1}^{k-1})^j~~,\] 
\[ U_3=n^{j+1}C_{n-1}^{k-1}C_{n-1}^k(C_{n-1}^k)^j/(C_n^k)^{j+1}~~,\] 
see (\ref{Lk}) and (\ref{Rk}). In this difference the product 
$U_2$ is irrelevant. Indeed, the highest power of $(n-1)$ multiplying 
any of the variables $q_{\nu}^{(0)}$ in $R_k$ is higher than the 
highest power of $(n-1)$ in $U_2$, see (\ref{Rk}). 

Set $(n-1)\cdots (n-j-1)=(n-1)^{j+1}+V$. 
Hence $U_1=((n-1)^{j+1}+V)C_{n-j-2}^{k-1}$. By Remark~\ref{divisible} 
the quantity $V$ is a polynomial divisible by $j(j+1)$. Therefore for $j=0$ 
one has $U_1=(n-1)C_{n-2}^{k-1}=nC_{n-1}^{k-1}C_{n-1}^k/C_n^k=U_3$, and 
for $j=-1$ one has $U_1=C_{n-1}^{k-1}=U_3$. 
Hence the CTs are divisible by 
$j(j+1)$. This completes the proof of  
Proposition~\ref{polynomial}. Now we settle Lemmas~\ref{A} and \ref{poldeg}. 

\medskip
\noindent
{\em Proof of Lemma~\ref{A}:}   
$1^0$. Using the equalities $C_{n-1}^{k-1}/C_n^k=k/n$ and 
$C_{n-1}^{k}/C_n^k=(n-k)/n$, one can present 
the equality $L_k=R_k$ (see (\ref{Lk}) and (\ref{Rk})) in the form 
\[ (n-1)\cdots (n-j-1)((-1)^jC_{n-j-2}^{k-1}+C_{n-j-2}^{k-2}q_{j-1}+\cdots 
+C_{n-j-2}^0q_{j-k+1})=\]
\begin{equation}\label{Lk=Rk}
=(n-k)C_{n-1}^{k-1}(k^j+\sum _{\nu =1}^{j-1}(n-k)^{\nu}k^{j-\nu}q_{\nu}+
(-1)^j(n-k)^j).
\end{equation}
Replace in (\ref{Lk=Rk}) the quantities $q_{\nu}$ by their expansions 
(\ref{qTaylor}). Consider the right-hand side $R_k$ 
of (\ref{Lk=Rk}) as 
a Laurent series in $\frac{1}{n-1}$. Observe that if the integer $k$ is 
bounded, then the following relations hold: 
\[ (n-k)C_{n-1}^{k-1}=\frac{(n-1)^k}{k!}+O((n-1)^{k-1})~~,~~
(n-k)^{\nu}k^{j-\nu}=k^{j-\nu}((n-1)^{\nu}+O((n-1)^{\nu -1}).
\]  

We use the last equality for $\nu =j-s$. The coefficient of $(n-1)^{j-s+k}$ in 
$R_k$ is of the form: 
\[ \frac{1}{k!}(k^sq_{j-s}^{(0)}+k^{s-1}q_{j-s+1}^{(1)}+\cdots 
+kq_{j-1}^{(s-1)}+{\mathcal H}+r),\] 
where ${\mathcal H}$ is a linear form in the variables $q_{\mu}^{(m)}$ with 
$\mu -m>j-s$ and $r$ is a real number. Hence the form ${\mathcal H}$ 
contains only  variables belonging to  
the union ${\mathcal J}_{j-s+1,s-2}$ (because $\mu \leq j-1$, $m\geq 0$) while 
the linear form $(1/k!)(k^sq_{j-s}^{(0)}+k^{s-1}q_{j-s+1}^{(1)}+\cdots 
+kq_{j-1}^{(s-1)})$ depends only on the variables in the set 
${\mathcal I}_{j-s,s-1}$. Hence $R_k$ depends only on the variables in the set 
${\mathcal J}_{j-s,s-1}$. 

$2^0$. Consider now the left-hand side $L_k$ of (\ref{Lk=Rk}). 
One can write 

\[ B(n,j):=(n-1)\cdots (n-j-1)=(n-1)^{j+1}\left( 1-\frac{e _1(j)}{n-1}+
\frac{e_2(j)}{(n-1)^2}-\cdots \right),\]   
see Notation~\ref{psikj}. For each $\nu =j-k+1,\ldots ,j$ the 
product $B(n,j)C_{n-j-2}^{\nu -j+k-1}$ is a  polynomial 
in the variable $(n-1)$ of degree $\nu +k$. More precisely,  

\begin{equation}\label{Bnj}
B(n,j)C_{n-j-2}^{\nu -j+k-1}=\frac{(n-1)^{\nu +k}}{(\nu -j+k-1)!}\left( 
1-\frac{e _1(\nu +k-1)}{n-1}+
\frac{e _2(\nu +k-1)}{(n-1)^2}-\cdots \right). 
\end{equation} 
Therefore, the coefficient of $(n-1)^{j-s+k}$ in the term 
$B(n,j)C_{n-j-2}^{\nu -j+k-1}q_{\nu}$ of 
$L_k$ is of the form 

\begin{equation}\label{formula} 
\frac{1}{(\nu -j+k-1)!}(q_{\nu}^{(\nu -j+s)}
-q_{\nu}^{(\nu -j+s-1)}e _1(\nu +k-1)+\cdots +
(-1)^{\nu -j+s}q_{\nu}^{(0)}e _{\nu -j+s}(\nu +k-1)).
\end{equation} 
The index $\nu$ takes the values $j-k+1,\ldots ,j-1$, see (\ref{Lk=Rk}). 
Hence $L_k$ is also a linear inhomogeneous form 
of the variables in the set ${\mathcal J}_{j-s,s-1}$.~~~~~$\Box$\\ 

\noindent
{\em Proof of Lemma~\ref{poldeg}:}
$1^0$. Consider equation $(A_{k,s-j-k})$. Recall that 
its unknown variables are 
the ones in the set ${\mathcal I}_{j-s,s-1}$. Present this equation in the 
form 
$\alpha _1q_{j-s}^{(0)}+\cdots +\alpha _sq_{j-1}^{(s-1)}=\beta +{\mathcal G}$ 
where the term $\beta$ depends on $j$ but not on the variables $q_{\nu}^{(i)}$ 
and ${\mathcal G}$ is a linear form in the variables $q_{\nu}^{(i)}$ from the 
set ${\mathcal J}_{j-s+1,s-2}$. 

$2^0$. The quantity $\beta$ 
is obtained  
by adding the terms $(n-1)^{j-s+k}$ from the Laurent series of the 
expressions : 
$A:=(n-1)\cdots (n-j-1)(-1)^jC_{n-j-2}^{k-1}$, $D:=-(n-k)C_{n-1}^{k-1}k^{j}$ 
and $W:=C_{n-1}^{k-1}(-1)^j(n-k)^{j+1}$ in equation (\ref{Lk=Rk}). 
The coefficient of  $(n-1)^{j-s+k}$ in $A$ equals 
$(-1)^{j+s}e _s(j+k-1)/(k-1)!$ (see equation (\ref{Bnj}) with $\nu =j$) 
which is a degree $2s$ polynomial in $j$, see Lemma~\ref{divisible}. 
Its coefficient in $W$ is a polynomial in $j$ of degree $s$. Indeed, 
\[ W=(-1)^jC_{n-1}^{k-1}((n-1)-(k-1))^{j+1}=
(-1)^jC_{n-1}^{k-1}\sum _{\gamma =0}^{j+1}C_{j+1}^{\gamma}(n-1)^{\gamma}
(k-1)^{j+1-\gamma},\] 
where $(-1)^jC_{n-1}^{k-1}$ is a  polynomial in $(n-1)$ of degree $k-1$ and 
$C_{j+1}^{j-s+1}=C_{j+1}^{s}$ is a polynomial in $j$ of degree $s$. 
As $k\leq s<j$, there is no term $(n-1)^{j-s+k}$ in $D$. Hence $\beta$ 
is a  polynomial in $j$ of degree $2s$.
 
$3^0$. The linear form ${\mathcal G}$ is obtained from certain expressions in 
both sides of equality (\ref{Lk=Rk}). First, 
one considers  the terms $B(n,j)C_{n-j-2}^{\nu -j+k-1}q_{\nu}$ in $L_k$ and 
their coefficients of $(n-1)^{j-s+k}$ given by formula (\ref{formula}). 
Recall that (see Lemma~\ref{A}) the index $\nu$ takes values $\geq j-s$. 
Set $\theta :=j-\nu$. Hence $\theta \leq s$. 
In formula (\ref{formula}) 
the term $\sigma :=q_{\nu}^{(\delta)}e _{\nu -j+s-\delta}=
q_{\nu}^{(\delta)}e _{s-\theta -\delta}$ 
is the product of the degree $2(s-\theta -\delta )$ 
polynomial $e_{s-\theta -\delta}$  in the variable $j$ 
(see Remark~\ref{divisible})  and of 
$\pm q_{j-\nu}^{(\delta)}=\pm q_{\theta}^{(\delta)}$ which is 
a polynomial of degree $2(\theta +\delta )$ by inductive assumption. Thus 
$\sigma$ (and, hence, the whole contribution 
of $L_k$ to the linear form ${\mathcal G}$)  is a   polynomial in $j$ 
of degree $2s$.

$4^0$. Secondly, consider in $R_k$ the product 

\[ S:=C_{n-1}^{k-1}(n-k)^{\nu +1}k^{j-\nu}q_{\nu}=
C_{n-1}^{k-1}k^{\theta}((n-1)+(1-k))^{\nu +1}
\sum _{\eta =0}^{\infty}q_{\nu}^{(\eta)}(n-1)^{-\eta}~.\] 

The quantity $q_{\nu}^{(\eta)}=\pm q_{\theta}^{(\eta)}$ is a  polynomial in 
$j$ of degree $2(\theta +\eta )$. 

\medskip 
Our goal now  is  to show that   
{\em the coefficient of $(n-1)^{k+\nu -r}$ in the product 
$C_{n-1}^{k-1}(n-k)^{\nu +1}$ is a  polynomial in $j$ of degree $r$}. 
(We prove this statement in $5^0$ below.) 
This implies that the coefficient of  
$(n-1)^{j-s+k}$ in $S$ is a finite sum of polynomials in $j$ of degrees
$\tau :=2(\theta +\eta )+r$. To obtain a term $(n-1)^{j-s+k}$ we multiply 
the terms $(n-1)^{k+\nu -r}$ and $(n-1)^{-\eta}$. In other words,  one has 
$(k+\nu -r)-\eta =j-s+k$, i.e. 
$\nu =j+r+\eta -s$ and $\tau =2(j-\nu +\eta )+r=2s-r<2s$. 
Thus the contribution of $R_k$ to 
the linear form ${\mathcal G}$ is a  polynomial in $j$ of degree $<2s$. The 
lemma is proved.

$5^0$. Proof of the latter statement.  
One has $(n-k)^{\nu +1}=
\sum _{r=0}^{\nu +1}C_{\nu +1}^{r}(1-k)^{r}(n-1)^{\nu +1-r}$ and 
$C_{\nu +1}^{r}$ is a degree $r$ polynomial in $\nu$, and therefore, also in 
the variable $j-\theta$ and thus  in the variable  $j$ as well. The binomial 
coefficient $C_{n-1}^{k-1}$ is a polynomial in $(n-1)$ of degree $(k-1)$ . 
Thus $C_{n-1}^{k-1}(n-k)^{\nu +1}$ is of the form 
$\sum _{r=0}^{k+\nu}d_{r}(n-1)^{k+\nu -r}$ where $d_{r}$ 
is a  polynomial in $j$ of degree $r$.~~~~~$\Box$ 

Now we finally return to our main results formulated in the introduction. 

\medskip
\subsection{Proof of Theorem~\protect\ref{MjNj}\protect\label{prMjNj}}

$1^0$. Consider the lower triangular matrix ${\mathcal M}$ 
whose $j$-th row contains 
the coefficients of the polynomial $M_j(x)$ (starting with the coefficient 
of the linear term) 
followed by zeros. Let us turn the Narayana triangle (\ref{NT}) into an 
infinite lower triangular matrix of the form 
\begin{equation}\label{NNT} 
\mathcal N=\begin{pmatrix}1&0&0&0&0&\cdots\\1&1&0&0&0&\cdots\\
1&3&1&0&0&\cdots\\
1&6&6&1&0&\cdots\\1&10&20&10&1&\cdots\\
\vdots&\vdots&\vdots&\vdots&\vdots&\ddots
\end{pmatrix}
\end{equation} 

Theorem~\protect\ref{MjNj} claims that the matrices 
${\mathcal M}$ and ${\mathcal N}$ coincide.  
Denote by ${\mathcal M}_l$ and ${\mathcal N}_l$ their $l$-th 
columns and by ${\mathcal M}^l$ and ${\mathcal N}^l$ their $l$-th diagonals 
(i.e. the sets of entries in positions $(r,r+l-1)$, $r=1,2,\ldots$ in 
${\mathcal M}$ and ${\mathcal N}$ respectively).

The polynomials $M_j(x)$ are monic and self-reciprocal  with positive 
coefficients by definition. 
Therefore, one has ${\mathcal M}_1={\mathcal N}_1$, 
${\mathcal M}^1={\mathcal N}^1$. 

$2^0$. Suppose that the first $m$ columns (and, hence, the first $m$ 
diagonals as well) 
of ${\mathcal M}$ coincide with the first $m$ columns (respectively, 
diagonals) of ${\mathcal N}$. The first $m$ entries of 
${\mathcal M}_{m+1}$ and of ${\mathcal N}_{m+1}$ vanish. 
Their next $m$ entries 
belong to the first $m$ diagonals, hence, they  coincide as well. 

$3^0$. The entries of ${\mathcal M}_{m+1}$ and the ones of 
${\mathcal N}_{m+1}$ (denoted respectively 
by ${\mathcal M}_{m+1,j}$ and ${\mathcal N}_{m+1,j}$) are the values  of 
polynomials $R_M^{m+1}$ and $R_N^{m+1}$ in $j$ of the same degree $2m$. 
For ${\mathcal M}_{m+1,j}$ this follows from 
Proposition~\ref{polynomial}, and for ${\mathcal N}_{m+1,j}$ 
it follows from the next formula for the Narayana numbers:

\begin{equation}\label{NN} 
N_{j,m+1}=\frac{1}{j}C_{j}^{m}C_{j}^{m+1}=
\frac{(j-1)(j-2)\cdots (j-m+1)}{m!}\frac{j(j-1)\cdots (j-m)}{(m+1)!}.
\end{equation} 

For $j=1,\ldots ,2m$ one has $R_M^{m+1}(j)=R_N^{m+1}(j)$, see $2^0$. 
For $j=0$ one has $N_{j,m+1}={\mathcal N}_{m+1,j}=0={\mathcal M}_{m+1,j}$. 
The first two equalities follow from formula (\ref{NN}), the last one can be 
deduced from part (2) of Proposition~\ref{polynomial}. This proposition 
implies that  ${\mathcal M}_{m+1,j}$ is divisible by $j(j-1)$ (recall that 
$M_j(x)=(-1)^{j-1}xQ_{j-1}^*(-x)$). 
Hence  $R_M^{m+1}(x)=R_N^{m+1}(x)$.~~~~~$\Box$



\subsection{Proof of Theorem~\protect\ref{interlace}
\protect\label{printerlace}}

$1^0$. For $n=2$ and $3$ all statements of the theorem can be checked 
directly. Observe that $N_n(x)/x$ has all  coefficients positive. By 
Corollary~\ref{hyperbol} all roots of $N_n(x)/x$ are real, hence negative, and 
$0$ is a simple root of $N_n(x)$. 

For any  even $n$ it follows from 
$N_n(x)/x$ being self-reciprocal, of odd degree 
and with positive coefficients that $N_n(-1)=0$. 
For $n$ odd the polynomial $N_n(x)$ does not vanish at $(-1)$.  
This can be proved by induction 
on $n$. Namely,  if $N_{n-2}(x)$ does not vanish at $(-1)$, then the same 
holds for $N_n(x)$, see (\ref{recur}).

$2^0$. We prove the rest of the theorem using induction on $n$. 
Suppose that its statements hold for  $N_2(x)$, $\ldots$, $N_{n-1}(x)$. 
Denote by $\xi _i$ the $i$-th root of $N_{n-1}(x)$, $\xi _i<\xi _{i+1}$, 
$\xi _{n-1}=0$. 
By part (3) of the theorem (proved for $N_{n-1}(x)$) 
one has $N_{n-2}(\xi _i)/\xi _i\neq 0$ 
and the sign of $N_{n-2}(\xi _i)/\xi _i$ changes alternatively. 
By equality (\ref{recur}) so does the sign of $N_n(\xi _i)/\xi _i$ as well 
(equality (\ref{recur}) implies that the signs of $N_{n-2}(\xi _i)/\xi _i$ 
and $N_n(\xi _i)/\xi _i$ are opposite). 
As $N_{n-2}(\xi _{n-2})/\xi _{n-2}>0$, one has 
$N_{n}(\xi _{n-2})/\xi _{n-2}<0$. The leading coefficient of 
the polynomial $N_{n}(x)/x$ is positive, 
hence it has a root in $(\xi _{n-2},0)$ and 
by self-reciprocity a root in the interval $(-\infty ,\xi _1)$ as well. 

The polynomial $N_n(x)/x$ is of degree $(n-1)$ and has 
at least one root in each of the intervals $(-\infty ,\xi _1)$, 
$(\xi _1,\xi _2)$, $\ldots$, $(\xi _{n-2},\xi _{n-1})$. 
Hence all these roots are simple 
(including $(-1)$ for $n$ even) and they interlace 
with the roots of $N_{n-1}(x)/x$. Thus the only root in common for 
$N_n(x)$ and $N_{n-1}(x)$ is $0$ which is  simple  for both of them, 
see part 1) of Remark~\ref{Narrem}.~~~~~$\Box$


\subsection{On root asymptotics of Narayana polynomials}
\label{asymptotics}
 
 In this subsection we prove Theorem~\ref{dens}.  Define 
$\Psi_n(x)=\frac{N_{n+1}(x)}{N_{n}(x)}$ and 
$\Theta_n(x)=\frac{N_n^\prime(x)}{nN_n(x)},\;n=1,2,\ldots$. 
Set $\Psi(x)=\lim_{n\to\infty}\Psi_n(x)$ and  
$\Theta(x)=\lim_{n\to\infty}\Theta_n(x)$ where these limits exist. $\Psi(x)$ 
is called the {\em asymptotic quotient} and $\Theta(x)$ is called the 
{\em asymptotic logarithmic derivative} of the sequence $\{N_n(x)\}$. 
 
\begin{lm}\label{ratio} The sequence $\{\Psi_n(x)\}$ of rational functions 
converges in $\mathbf C\setminus \mathbf R_{\le 0}$ to the function 
$\Psi(x)=x+1+2\sqrt{x}=(\sqrt{x}+1)^2$ where $\sqrt{x}$ is the usual branch of 
the square root which attains positive values for positive $x$. (Here $\mathbf 
R_{\le 0}$ denotes the  half-axis of non-positive numbers.) 
 \end{lm}
 
 \noindent
{\em Proof of Lemma~\ref{ratio}:} We need to invoke the classical result of 
H.~Poincar\'e, see \cite[p. 287 - 298]{Ge} and Theorem~\ref{th:Poin} 
and Remark~\ref{HPrem} below.  
Indeed, dividing the recurrence relation (\ref{recur}) by $n+1$ 
we obtain the normalized reccurence
\begin{equation}\label{recurN} 
N_n(x)-\frac{2n-1}{n+1}(x+1)N_{n-1}(x)+\frac{n-2}{n+1}(x-1)^2N_{n-2}(x)=0. 
\end{equation} 
Taking  limits of its coefficients when $n\to \infty$ we get from 
Poincar\'e's theorem that the asymptotic  quotient  $\Psi(x)$ for  
each $x$ when $\Psi(x)$ exists should satisfy the following quadratic 
characteristic equation: 
\begin{equation}\label{eq:limrec}
\Psi^2(x)-2(x+1)\Psi+(x-1)^2=0.
\end{equation}Ê 
The exceptional set $E\subset  \mathbf C$  (called {\em the equimodular 
discriminant}, see \cite{Bi}) of those values of $x$ for which the equation 
(\ref{eq:limrec}) does not hold consists of all  $x$ for which two different 
solutions $\Psi_1(x)$ and $\Psi_2(x)$ of (\ref{eq:limrec}) have the same 
absolute value. 

Let us show now that in the considered case $E=\mathbf R_{\le 0}$. Indeed, 
two solutions of (\ref{eq:limrec}) are given by $\Psi_1(x)=x+1+2\sqrt{x}$ and 
$\Psi_2(x)=x+1-2\sqrt{x}$ for some choice of the branch of square root. One 
can easily check that if $|\Psi_1(x)|=|\Psi_2(x)|$ for some value of $x$ then 
$x+1$ is orthogonal to $\sqrt{x}$ as two vectors in $\mathbf R^2\simeq \mathbf 
C$. Denoting $\sqrt{x}=A+iB$ one gets $x+1=A^2-B^2+1+i(2AB)$ and the latter 
orthogonality condition is given by the relation: 
$$A(A^2-B^2+1)+B(2AB)=0\Longleftrightarrow A(A^2+B^2+1)=0\Longleftrightarrow 
A=0.$$ 
But  the real part $A$ of $\sqrt{x}$ vanishes if and only if $x$ is a 
negative real number implying  $E=\mathbf R_{\le 0}$.  

Thus  the asymptotic quotient $\Psi(x)$ satisfies the equation 
(\ref{eq:limrec}) in  $\mathbf C\setminus\mathbf R_{\le 0}$. To show that  
$\Psi(x)=x+1+2\sqrt{x}$, i.e. it chooses the correct branch of solutions to  
(\ref{eq:limrec}) we check that $\Psi(1)=1+1+2=4$ and we prove that $\Psi(x)$ 
is continuous in $\mathbf C\setminus \mathbf R_{\le 0}$. 

Indeed, one knows that $N_n(1)=\sum_{k=1}^nN_{n,k}=Cat_n$ where 
$Cat_n=\frac{1}{n+1}\binom{2n}{n}$ is the $n$-th Catalan number, see 
\cite {Su}.  Thus $\Psi_n(1)=\frac{N_{n+1}(1)}{N_{n}(1)}=
\frac{2(n+1)(2n+1)}{(n+2)(n+1)}.$ Therefore, $\Psi(1)=\lim_{n\to \infty}
\Psi_n(1)=4$. In order to prove the required continuity (and analiticity) we 
use the well-known Montel's theorem claiming that a locally bounded family of 
analytic functions contains a subsequence converging to  an analytic function, 
see e.g \cite{Da}. For us it is technically easier to work with the sequence 
of inverses to $\Psi_n(x)$, i.e.  with the sequence $\{\frac{N_n(x)}
{N_{n+1}(x)}\}$. We show that $\{\frac{N_n(x)}{N_{n+1}(x)}\}$ is locally 
bounded in any compact domain separated from $\mathbf R_{\le 0}$. Thus since 
by Poincar\'e's theorem $\{\frac{N_n(x)}{N_{n+1}(x)}\}$  is pointwise 
converging in $\mathbf C\setminus \mathbf R_{\le 0}$ one gets  by Montel's 
theorem that $\{\frac{N_n(x)}{N_{n+1}(x)}\}$ converges to an analytic function 
implying the same fact for the sequence $\{\Psi_n(x)\}$. Indeed, 
by (2) and (3) of Theorem~\ref{interlace} all roots of each $N_n(x)/x$ are 
simple, strictly negative and they interlace with that of $N_{n+1}(x)/x$. 
Therefore,  the  partial fractional decomposition of $\frac{N_n(x)}
{N_{n+1}(x)}$ has the form 
$$\frac{N_n(x)}{N_{n+1}(x)}=\sum_{i=1}^{n}\frac{\gamma_i}{x-a_{i,n+1}},$$
where every $\gamma_i$ is positive with $\sum_{i=1}^n\gamma_i=1$ and 
$\{a_{1,n+1},\ldots ,a_{n,n+1}\}$ 
is the set of all non-vanishing (and therefore 
strictly negative) roots of $N_{n+1}(x)$. If $x\in \mathbf C$ is any point 
denote by $\nu(x)$ its  Eucledian distance to $\mathbf R_{\le 0}$. Let us 
check that for any 
$x\in \mathbf C\setminus \mathbf R_{\le 0}$ and  for any positive integer 
$n$ one has 
$$\left| \frac{N_n(x)}{N_{n+1}(x)} \right|\le \frac {1}{\nu(x)},$$
which immediately implies the required local boundedness. Indeed, 
$$\left|\frac{N_n(x)}{N_{n+1}(x)}\right|=\left|\sum_{i=1}^{n}\frac{\gamma_i}
{x-a_{i,n+1}}\right|\le \sum_{i=1}^{n}\frac{\gamma_i}{|x-a_{i,n+1}|}\le 
\sum_{i=1}^{n}\frac{\gamma_i}{\nu(x)}=\frac{1}{\nu(x)}.$$
~~~~~$\Box$

 \begin{lm}\label{quotient} The sequence $\{\Theta_n(x)\}$ of rational 
functions converges in $\mathbf C\setminus \mathbf R_{\le 0}$ to the function 
$\Theta(x)=\frac{1}{x+\sqrt{x}}$. 
 \end{lm}
 
 \noindent 
{\em Proof of Lemma~\ref{quotient}:} 
Indeed, it is known, see \cite{ST}, that in case when both the asymptotic 
quotient $\Psi(x)$ and the asymptotic ratio $\Theta(x)$ exist and have 
continuous first derivatives in some open domain of $\mathbf C$ they  satisfy 
there the relation 
\begin{equation}\label{rel}
\Theta(x)=\frac{\Psi'(x)}{\Psi(x)}.
\end{equation}
A short sketch of its proof is as follows. Consider the difference 
$$\frac{\Psi^\prime_n(x)}{\Psi_n(x)}- \Theta_n(x)=\frac{p_{n+1}^\prime(x)}
{p_{n+1}(x)}-\left(1+\frac{1}{n}\right)\frac{p_{n}^\prime(x)}{p_{n}(x)}.$$
Then one has, 
$$\frac{\Psi^\prime(x)}{\Psi(x)}-\Theta(x)=\lim_{n\to\infty} 
\left(\frac{\Psi^\prime_n(x)}{\Psi_n(x)}- \Theta_n(x)\right)=
\lim_{n\to\infty}(n+1)\left(\frac{p_{n+1}^\prime(x)}{(n+1)p_{n+1}(x)}-
\frac{p_{n}^\prime(x)}{np_{n}(x)}\right)=0.$$

Applying (\ref{rel})  to the formula for $\Psi(x)$ in Lemma~\ref{ratio} we get 
$$\Theta(x)=\frac{\Psi'(x)}{\Psi(x)}=\frac{2(\sqrt{x}+1)}
{2\sqrt{x}(\sqrt{x}+1)^2}=\frac{1}{\sqrt{x}(\sqrt{x}+1)}=
\frac{1}{x+\sqrt{x}}.$$~~~~~$\Box$ 

\begin{nota} Given a finite measure $\mu$ supported on $\mathbf C$ define 
its Cauchy transform $\C_\mu(x)$ as 
$$\C_\mu(x)=\int_{\mathbf C}\frac{d\mu(\zeta)}{x-\zeta}.$$
\end{nota}

The Cauchy transform of the measure is analytic outside its support, 
its $\frac{\partial}{\partial \bar z}$-derivative coincides with the original 
measure and it has many more important properties, see e.g. \cite {Ga}. 
Notice that for any polynomial $P(x)$ of degree $l$ the Cauchy transform 
$\C_{\mu_P}$Êof its root-counting measure $\mu_P$ is given by the formula 
$\C_{\mu_P}(x)=\frac{P^\prime(x)}{l\cdot P(x)}$. Therefore, given a polynomial 
sequence $\{P_n(x)\}, \deg P_n(x)=n, n=1,2,\ldots$ we have that the Cauchy 
transform $\C_\mu$ of its asymptotic root-counting measure 
$\mu=\lim_{n\to\infty}\mu_{P_n}$ (if the latter measure exists) coincides 
with the limit $\C_\mu(x)=\lim_{n\to\infty}\frac{P_n^\prime(x)}{nP_n(x)}$.

The last result which we need to settle Theorem~\ref{dens} and which is a 
particular case of Theorem 3.1.9 of \cite {Ho} is as follows.

\begin{lm}\label{density} If $\rho_\mu(x), x\in \mathbf R_{\le 0}$ is the 
density of a  finite measure $\mu$   supported on $\mathbf R_{\le 0}$ and 
$\C_\mu(x), x \in \mathbf C \setminus \mathbf R_{\le 0}$ is its Cauchy 
transform then for any $x\in \mathbf R_{\le 0}$ one has
$$\rho_\mu(x)=\frac{i}{2\pi} \left(\overline{\C_\mu^+(x)-\C_\mu^-(x)}
\right), $$
where $\C_\mu^+(x)=\lim_{t\to x}\C_\mu(x)$ with $t$ belonging to the upper 
halfplane,  $\C_\mu^-(x)=\lim_{t\to x}\C_\mu(x)$ with $t$ belonging to the 
lower halfplane, and $\overline Z$ denotes the usual conjugate of $Z$.
\end{lm} 

Finally, applying the latter formula to our case we get the required density 
formula in the statement of Theorem~\ref{dens}. Namely, 
$$\rho(x)=\frac{i}{2\pi}\left(\overline{\frac{1}{x+i\sqrt{-x}}-
\frac{1}{x-i\sqrt{-x}}}\right)
=\frac{i}{2\pi(x^2-x)} \left(\overline{(x-i\sqrt{-x})-(x+i\sqrt{-x})}\right)=$$
$$=\frac{i}{2\pi(x^2-x)}\cdot\overline{-2i\sqrt{-x}}=\frac{-i^2}{\pi}
\frac{\sqrt{-x}}{x^2-x}=\frac{1}{\pi(1-x)\sqrt{-x}}.$$
Integrating the obtained formula for $\rho(x)$ one gets the expression for 
the distribution function $\kappa(x)$ from the statement of 
Theorem~\ref{dens}, see Fig.~1.~~~~~$\Box$

\section {Appendix. Poincar\'e's theorem, \cite[p. 287]{Ge}}\label{Poincare}

{\bf 1. Set-up.} We consider a linear homogeneous difference equation of
order $k$
\begin{equation}\label{eq:1}
f(t+k)+a_{1}f(t+k-1)+a_{2}f(t+k-2)+\cdots+a_{k}f(t)=0
\end{equation}
with constant coefficients. Denote by
$$\la_{1},\ldots,\la_{k}$$
the roots of the characteristic equation
$$\la^k+a_{1}\la^{k-1}+a_{2}\la^{k-2}+\cdots+a_{k}=0,$$
  and assume that $\la_{i}$ have different absolute values denoted by 
  $\xi_{i}=\vert \la_{i}\vert.$  We can then assume that
  \begin{equation}\label{eq:2}
  \xi_{1}>\xi_{2}>\cdots>\xi_{k}.
  \end{equation}
  Since the roots $\la_{i}$ are all distinct then the general solution
  of (\ref{eq:1}) is given by
  $$f(t)=C_{1}\la_{1}^t+C_{2}\la_{2}^t+\cdots+C_{k}\la_{k}^t.$$
Let us  choose one single solution of (\ref{eq:1}), i.e. assign some
  fixed values to the constants $C_{i}$. Let $C_{p}$ be the first
  nonvanishing  among $C_{i}$'s, i.e.
  $$C_{1}=C_{2}=\cdots=C_{p-1}=0\quad{\text {and}}\quad C_{p}\neq 0.$$

  Then one can show that
  $$\lim_{t\to\infty}\frac{f(t+1)}{f(t)}=\la_{p}$$
  where $f(t)$ is the considered solution of (\ref{eq:1}). Indeed, one has
   $$ \frac {f(t+1)}{f(t)}=\frac
   {C_{p}\la_{p}^{t+1}+C_{p+1}\la_{p+1}^{t+1}+\cdots+C_{k}\la_{k}^{t+1}}
   {C_{p}\la_{p}^{t}+C_{p+1}\la_{p+1}^{t}+\cdots+C_{k}\la_{k}^{t}}
   $$
   For $t\to\infty$ one has through (\ref{eq:2}) that the limits of
   $$\left(\frac {\la_{p+1}}{\la_{p}}\right)^t, \left(\frac
   {\la_{p+2}}{\la_{p}}\right)^t, \ldots , \left(\frac
   {\la_{k}}{\la_{p}}\right)^t$$
   equal $0$, and we obtain
    $$\lim_{t\to\infty}\frac{f(t+1)}{f(t)}=\la_{p}$$
    that is the following statement is valid

\begin{tm} If $f(t)$ is an arbitrary (nontrivial) solution of the
equation (\ref{eq:1}) then the limit of the ratio $\frac {f(t+1)}{f(t)}$ for
$t\to\infty$ equals one the roots of the characteristic equation if
only these roots have distinct absolute values.

     \end{tm}

    We cannot say which root will be involved without the knowledge of
    the solution. One can only say that this root will have the number $p$
    of the first nonvanishing coefficient $C_{p}$ in the considered
    solution. The following theorem of Poincar\'e presents a
    generalization of this statement.

\medskip 
   {\bf  2. Poincar\'e's theorem.}

\begin{tm}\label{th:Poin} If the coefficients $P_{i}(t),\;i=0,1,\ldots,k-1$ 
of a linear homogeneous difference equation
\begin{equation}\label{eq:3}
f(t+k)+P_{k-1}(t)f(t+k-1)+P_{k-2}(t)f(t+k-2)+\cdots+P_{0}(t)f(t)=0
\end{equation}
have limits $\lim_{t\to\infty}P_{i}(t)=a_{i},\;\;i=0,1,\ldots,k-1$
and if the roots of the characteristic equation
\begin{equation}\label{eq:4}
\la^k+a_{k-1}\la^{k-1}+\cdots+a_{0}=0
\end{equation}
have different absolute values then the limit of the ratio
$\frac{f(t+1)}{f(t)}$ for $t\to \infty$ of any solution $f(t)$ of the
equation (\ref{eq:3}) equals one of the roots
$$\la_{1},\la_{2},\ldots,\la_{k}$$
of the equation (\ref{eq:4}), i.e.
$$\lim_{t\to\infty}\frac{f(t+1)}{f(t)}=\la_{p}.$$

     \end{tm}

\begin{rem}\label{HPrem}
In the present paper the role of the quantities $f(t+i)$ is played by the 
values of the polynomials $N_n$ for each $x$ fixed. The parameter $t$ 
is discrete -- this is the index $n$.
\end{rem}

\end{document}